\documentclass[final,1p,times]{elsarticle}

\usepackage{amssymb}
\usepackage{amsmath}

\usepackage{caption}
\usepackage{subfig}
\usepackage{bm}
\usepackage{booktabs}
\usepackage[numbers]{natbib}
\usepackage{float}
\usepackage{hyperref}
\journal{}

\begin{document}
	
	\begin{frontmatter}
		
		\title{Learning high-accuracy numerical schemes for hyperbolic equations on coarse meshes}
		
		\author[a]{Jinrui Zhou}
		\ead{zhou1218@std.uestc.edu.cn}
		
		\author[a]{Yiqi Gu}
		\ead{yiqigu@uestc.edu.cn}
		
		\author[a]{Hua Shen}
		\ead{huashen@uestc.edu.cn}
		
		\author[a]{Liwei Xu\corref{cor1}}
		\ead{xul@uestc.edu.cn}
		
		\author[b]{Juan Zhang}
		\ead{zhang_juan@buaa.edu.cn}
		
		\author[c]{Guanyu Zhou}
		\ead{zhoug@uestc.edu.cn}

		\affiliation[a]{
			organization={School of Mathematical Sciences, University of Electronic Science and Technology of China},
			city={Chengdu},
			state={Sichuan},
			postcode={611731},
			country={China}
		}
		\affiliation[b]{
			organization={Institute of Artificial Intelligence, Beihang University},
			city={Beijing},
			postcode={100191},
			country={China}
		}
		\affiliation[c]{
			organization={Institute of Fundamental and Frontier Sciences, University of Electronic Science and Technology of China},
			city={Chengdu},
			state={Sichuan},
			postcode={611731},
			country={China}
		}
		\cortext[cor1]{Corresponding author}
		
		\begin{abstract}
			When solving partial differential equations using classical schemes such as finite difference or finite volume methods, sufficiently fine meshes and carefully designed schemes are required to achieve high-order accuracy of numerical solutions, leading to a significant increase in computational costs, especially for three-dimensional (3D) time-dependent problems. Recently, machine learning-assisted numerical methods have been proposed to enhance accuracy or efficiency. In this paper, we propose a data-driven finite difference numerical method to solve the hyperbolic equations with smooth solutions on coarse grids, which can achieve higher accuracy than classical numerical schemes based on the same mesh size. In addition, the data-driven schemes have better spectrum properties than the classical schemes, although the spectrum properties are not explicitly optimized during the training process. Numerical examples are presented to demonstrate the accuracy and efficiency of the proposed method, as well as its good performance on dispersion and dissipation. 
		\end{abstract}

		\begin{keyword}
			Neural network \sep Data-driven \sep Finite difference scheme \sep Euler equations
			
		\end{keyword}
		
	\end{frontmatter}

		\section{Introduction}
		\label{section1}
		Fluid dynamics is a fundamental discipline that serves a wide range of scientific and engineering fields \cite{ref00, ref01}. The practical flows of interest typically have multiple scales in both space and time. If we aim to obtain a numerical solution encompassing a full range of spatial and temporal scales, the computational costs become expensive, and therefore the development of high-accuracy numerical solvers seems to be desirable. The central difference scheme, upwind difference scheme \cite{ref1}, modified upwind difference scheme optimized for dispersion and dissipation \cite{ref2, ref3, ref6} and Gaussian radial basis function method \cite{ref7} are commonly employed to address problems without shock. Nevertheless, the discretization errors of these schemes escalate as the mesh becomes coarser.
		
		In recent years, the neural network method for solving partial differential equations has attracted much attention due to its powerful approximation ability \cite{ref51, ref52}. The Physics-Informed Neural Networks (PINNs) \cite{ref53}, Deep Ritz method (DRM) based on the variational principle \cite{ref54}, Least-Squares ReLU Neural Network (LSNN) \cite{ref55}, Structure-preserving neural networks with asymptotic-preserving (AP) loss functions\cite{ref531, ref532}, and various other methods \cite{ref56, ref57} have been proposed. Although these methods exhibit good performance for high-dimensional problems, their accuracy has no advantages compared to the classical method, and the repeated training process is required for different problems. Therefore, the hybrid numerical method, by replacing some parts of numerical solvers with neural networks to obtain a balance between computational cost and precision, has also been developed in \cite{ref13, ref12, ref15, ref8, ref81, ref14}. Bezgin et al. \cite{ref10} and Xue et al. \cite{ref11} improved the WENO schemes with fully connected neural networks by learning the weight function and the smoothness detector, respectively. Manrique de Lara and Ferrer \cite{ref9, ref91} have used neural networks to accelerate the high-order discretization Galerkin method by introducing a corrective forcing to the low-order solution. The data-driven discretization schemes learned by convolutional neural networks (CNN) have been proposed in \cite{ref13}, which are optimized to satisfy the equations on a low-resolution grid. This method has been applied to the passive scalar advection in \cite{ref14} and incompressible Navier-Stokes equations in \cite{ref15}.
		
		In this paper, we develop a weight-learning neural network(WLNN) to generate finite difference coefficients during the reconstruction of the interface flux with high accuracy. Inspired by the design strategy of traditional numerical schemes, we adjust the input of the network. Instead of using global data on all grid points \cite{ref13}, we use several local points around the interface flux as the input of the network, which is in much lower dimensions and can be trained more efficiently. We consider as an illustration a six-point scheme widely utilized for practical applications. Under the consistency condition, the coefficients have two constraints, which are strictly enforced in hard constraint layers. The data-driven method is trained offline and is applied to solve various scalar hyperbolic problems and the Euler equations in both two and three dimensions. These results indicate that the proposed method exhibits good spectral properties and outperforms the classical schemes, especially on coarse meshes. In addition, the proposed method is more efficient than the classical scheme for the Taylor-Green vortex problem under the same accuracy.

		The remainder of this paper is organized as follows. In Section \ref{section2}, we review the usual preliminaries to reconstruct interface flux. In Section \ref{section3}, the weight-learning neural network and the training procedure are described. Numerical examples are provided in Section \ref{section4}. Finally, we conclude this paper in Section \ref{section5}.
		
		\section{Preliminaries}
		\label{section2}
		Reconstructing the numerical flux is an important procedure during the design of the finite difference schemes solving the hyperbolic equations. For instance, with the six-point scheme, we describe how weights are constrained by the consistency conditions. This approach can be applied to other schemes in a similar process, such as the four-point scheme.
		\subsection{Scalar hyperbolic equation}
		\label{sec:2.1}
\begin{figure}
			\centering
			\includegraphics[scale=0.5]{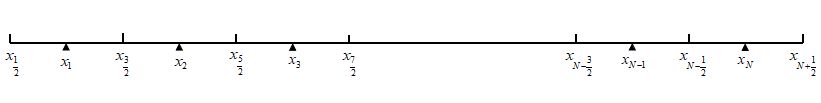}
			\caption{The discretization of the spatial domain.}
			\label{fig:2.1}
		\end{figure}

		Given flux $f$, we intend to find $u(x, t):[0, 1] \times [0,+\infty)\rightarrow\mathbb{R}$ such that
		\begin{align}
			\frac{\partial u}{\partial t} + \frac{\partial f(u)}{ \partial x} = 0, \quad x\in[0,1],~t \in [0, +\infty).
			\label{eq:2.1}
		\end{align}
		Let $N_x$ be a positive integer. We divide the interval $[0,1]$ into $N_x$ cells with the cell size $\Delta x = \frac{1}{N_x}$. For $j=1, \dots, N_x$, we denote by $I_j = [x_{j - \frac{1}{2}},x_{j + \frac{1}{2}}]$ the cell and by $x_j = (j-\frac{1}{2})\Delta x$ the cell center (see Fig. \ref{fig:2.1}). The conservative finite difference scheme can be derived as
		\begin{align}
			\frac{d u_i(t)}{d t} = -  \frac{\hat{f}_{i + \frac{1}{2}} - \hat{f}_{i - \frac{1}{2}}}{\Delta x},
			\label{eq:2.2}
		\end{align}
		where $u_i(t) := u(x_i, t)$ and the numerical flux $\hat{f}_{i\pm\frac{1}{2}}$ is an approximation to the flux function $h(x)$ implicitly defined as $f(x) = \frac{1}{\Delta x} \int_{x - \frac{\Delta x}{2} }^{x +  \frac{\Delta x}{2}} h(\xi) d \xi$ at $x_{i\pm\frac{1}{2}}$.
		Let $H(x) = \int_{0}^{x} h(\xi) d\xi$ be the primitive function of $h(x)$. Then $H_{i + \frac{1}{2}} := H(x_{i + \frac{1}{2}})$ can be expressed by
		\begin{align}
			H_{i + \frac{1}{2}} = \sum_{j=1}^{i} \int_{x_{j - \frac{1}{2}}}^{x_{j + \frac{1}{2}}} h(\xi) d \xi = \Delta x \sum_{j=1}^{i} f_j,
			\label{eq:2.3}
		\end{align}
		where $f_j:=f(u(x_j, t))$. Therefore, $h(x_{i + \frac{1}{2}})$ can be approximated as
		\begin{equation}
			\begin{aligned}
				h(x_{i + \frac{1}{2}}) = \frac{\partial H}{\partial x} \bigg|_{x=x_{i+\frac{1}{2}}}  = \sum_{l=1}^{7} a_l H_{i + l-\frac{7}{2}} + O(\Delta x ^k),
				\label{eq:2.4}
			\end{aligned}
		\end{equation}
		where $a_l$ is the coefficient of $H_{i + l-\frac{7}{2}}$. Substituting Eq. (\ref{eq:2.3}) to the Eq. (\ref{eq:2.4}), we obtain the reconstruction formula of $h(x_{i + \frac{1}{2}})$:
		\begin{equation}
			\begin{aligned}
				h(x_{i + \frac{1}{2}}) =  \sum_{l = 1}^{6} w_l f_{i-3+l} + O(\Delta x^k).
				\label{eq:2.5}
			\end{aligned}
		\end{equation}
		Note that $w_l = -\sum_{m = 1}^{l}a_l \Delta x$.
		To satisfy the consistency condition, we require $k \ge 2$ in Eq. (\ref{eq:2.5}) and the following the constraints of weights $w_l$: 
		\begin{subequations}
			\begin{align}
				w_1 + w_2 + w_3 + w_4 + w_5 + w_6 & = 1,	\label{eq:2.6a}		\\
				-5w_1 - 3w_2 - w_3 + w_4 + 3w_5 + 5w_6 & = 0.
				\label{eq:2.6b}
			\end{align}
		\end{subequations}	
		Hence, we obtain the numerical flux $\hat{f}_{i + \frac{1}{2}} := \sum_{l = 1}^{6} w_l f_{i-3+l}$ with the weights satisfying Eq. (\ref{eq:2.6a}) and (\ref{eq:2.6b}).
		
		There exist many classical schemes, such as 6th-order central scheme $\bm{w} = [\frac{1}{60}$, $-\frac{8}{60}$, $\frac{37}{60}, \frac{37}{60}$, $-\frac{8}{60}$, $\frac{1}{60}]$, 5th-order upwind scheme $\bm{w}=[\frac{1}{30}$, $ -\frac{13}{60}$, $ \frac{47}{60}$, $ \frac{9}{20}$, $ -\frac{1}{20}, 0]$, etc. Although these schemes have high-order convergence rates, the numerical accuracy or spectral properties are not optimal for complex problems on coarse grids. This study seeks to utilize a neural network to generate "optimal" weights for various problems, potentially outperforming classical schemes on coarse meshes. The neural network is trained using a high-accuracy dataset.
		
		\subsection{Hyperbolic conservation system}
		\label{sec:2.2}
		Let $\bm{U}=[U^0, U^1, U^2]$ $:[0, 1] \times [0, +\infty) \rightarrow \mathbb{R}^3$ be a vector-valued function. Given the vector flux function $\bm{F}(\bm{U})=[F^0(\bm{U}), F^1(\bm{U}), F^2(\bm{U})]$, we consider the 1D hyperbolic conservation system:  
		\begin{align}
			\frac{ \partial \bm{U}}{\partial t} + \frac{\partial \bm{F}}{ \partial x} = 0, \quad x\in[0, 1], t\in[0, +\infty).
			\label{eq:2.2.1}
		\end{align}

		The spatial domain is discretized uniformly, which is the same as that of Section \ref{sec:2.1}.  Then the conservative finite difference scheme for Eq. (\ref{eq:2.2.1}) reads: 
		\begin{align}
			\frac{ d \bm{U}_{i}(t)}{d t} =
			\frac{ \hat{\bm{F}}_{i+\frac{1}{2}} - \hat{\bm{F}}_{i-\frac{1}{2}}}{ \Delta x}, 
			\label{eq:2.2.2}
		\end{align}
		where $\bm{U}_{i}(t) := \bm{U}(x_i, t)$ and $\hat{\bm{F}}_{i\pm \frac{1}{2}}=[\hat{F}_{i\pm \frac{1}{2}}^0, \hat{F}_{i\pm \frac{1}{2}}^1, \hat{F}_{i\pm \frac{1}{2}}^2]$ is the numerical flux at the cell boundary. Similar to the reconstruction process of scalar equations, $\hat{F}_{i+\frac{1}{2}}^m$ is approximated by
		\begin{equation}
			\hat{F}_{i+\frac{1}{2}}^m = \sum_{l=1}^{6} w_l^m F_{i - 3 + l}^m \quad  (m = 0, 1, 2), 
			\label{eq:2.2.3}
		\end{equation}
		where the weights $\{w^m_l\}_{l=1}^{6}$ ($m=0, 1, 2$) satisfy Eq. (\ref{eq:2.6a}) and (\ref{eq:2.6b}).
		
		\section{Learning high-accuracy numerical schemes}
		\label{section3}
		In this section, we first introduce the architecture of neural networks. Then, the process of reconstructing the interface flux by the WLNN method is described. After that, we formulate the training dataset and the loss function.
		
		\subsection{The architecture of neural networks}
		\label{section3.1}
		We introduce the notations for neural networks, which are a family of nonlinear compound functions with a set of parameters that maps the input variables $\bm{x}$ to the output $\bm{y}=f(\bm{x}; \theta)$. The input and output variables can be either scalar or vector. A fully-connected neural network consists of three parts: input layer, hidden layers, and output layer:
		\begin{equation}
			\begin{aligned}
				& \text{Input layer} &
				& \bm{x}^0 = \bm{x},    \\
				& \text{Hidden layer} &
				& \bm{x}^s = \sigma \left( \bm{W}^{s} \bm{x}^{s-1} + \bm{b}^{s} \right), \quad s=1, \cdots, L-1   \\
				& \text{Output layer} &
				& \bm{y} = \sigma \left( \bm{W}^{L}\bm{x}^{L-1} + \bm{b}^{L} \right),
				\notag
				\label{eq:3.1.1}
			\end{aligned}
		\end{equation}
		where $\bm{W}^{s}$ and $\bm{b}^{s}$ are weight matrix and bias vector, respectively, linking the layers of $(s-1)$ and those of $s$. Neurons between adjacent layers are densely connected. The activation function $\sigma(\cdot)$ applied on a vector acts on every element of the vector. In general, the activation function in the last layer uses the identity function. We need to design a loss function $\mathcal{L}(\bm{y}, \hat{\bm{y}})$, which measures the difference between the output $\bm{y}$ and target $\hat{\bm{y}}$. With a proper optimization algorithm, the neural network gradually learns the mapping relation between $\bm{x}$ and $\hat{\bm{y}}$.
		
		\subsection{The weight-learning neural network}
		\label{section3.2}
		We introduce the flux splitting method and the architecture of the neural network.
		To construct a robust scheme, the flux of Eq. (\ref{eq:2.1}) is split into two parts \cite{ref16}
		\[
		f(u) = f^+(u) + f^-(u), 
		\]
		where $\frac{d f^+(u)}{d u} \ge 0$ and $\frac{d f^-(u)}{d u} < 0$. For instance, we can use the local Lax-Friedrichs (LLF) flux splitting method, i.e., 
		\begin{align}
			f^{\pm}(u) = \frac{1}{2}(f(u) \pm \alpha u), \notag
		\end{align}
		where $\alpha = \max \left| \frac{d f(u)}{ d u} \right|$ (the maximum is taken over several local points near the interface). Hence the interface flux $\hat{f}_{i+\frac{1}{2}}$ is decomposed as follows:  
		\begin{align}
			\hat{f}_{i+\frac{1}{2}} = \hat{f}^+_{i+\frac{1}{2}} + \hat{f}^-_{i-\frac{1}{2}},	\notag
		\end{align}
		where $\hat{f}^+_{i+\frac{1}{2}}$ and $\hat{f}^-_{i+\frac{1}{2}}$ represent the numerical flux of positive and negative values at $x_{i+\frac{1}{2}}$, respectively.
		
		Taking the six-point scheme as an example, we aim to learn a neural network that can yield a good six-point interpolation of $\hat{f}^+_{i+\frac{1}{2}}$ (resp. $\hat{f}^-_{i+\frac{1}{2}}$). More precisely, given the input $\{f^+_{i-2}, \cdots, f^+_{i+3}\}$, the neural network will output a group of weights $\bm{w} = \{w_{1}, \cdots, w_{6}\}$ such that $\sum_{l = 1}^{6} w_l f^+_{i - 3 + l}$ is the desired approximation of $\hat{f}^+_{i+\frac{1}{2}}$. 
		
		For simplicity, let us drop the superscript $+$. We denote by $\phi$, a $\mathbb{R}^6 \to \mathbb{R}^6$ mapping, the WLNN, i.e., 
		\begin{align}
			\bm{w}^{NN} = \phi(f_{i-2}, f_{i-1}, f_{i}, f_{i+1}, f_{i+2}, f_{i+3} ; \theta),
			\label{eq:3.2.2}
		\end{align}
		where the 6-point stencil $\{f_{i-2}, f_{i-1}, f_{i},$ $ f_{i+1}, f_{i+2}, f_{i+3}\}$ is the input, the group of weights $\bm{w}^{NN}=[w_1^{NN}, w_2^{NN}, \dots, w_6^{NN}]$ is the output, and $\theta$ represents the parameters of the neural network. 
		
\begin{figure}
			\centering
			\includegraphics[scale=0.35]{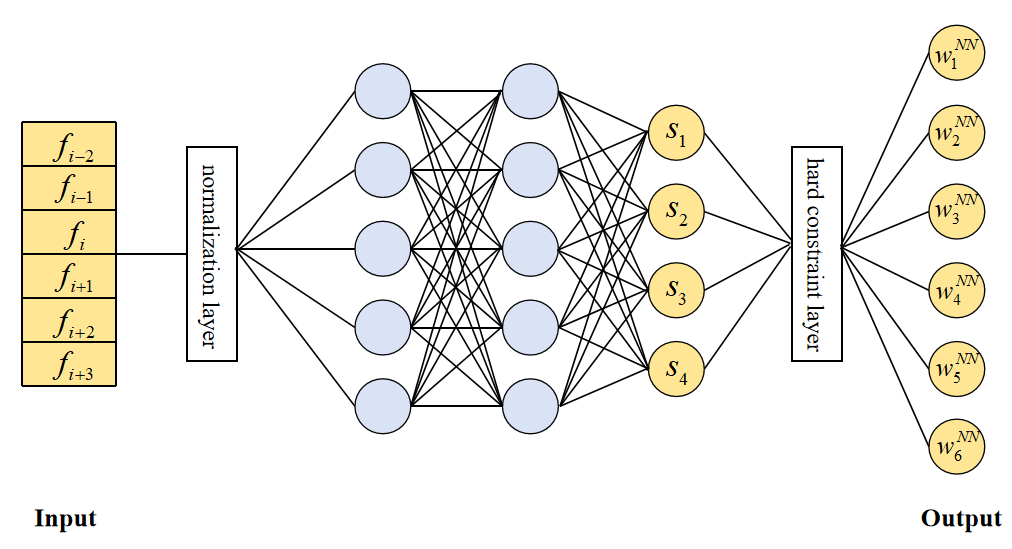}
			\caption{Schematic of the neural network architecture.}
			\label{fig:3.1}
		\end{figure}

		The architecture of WLNN is described below. 
		First, we pass the input to a normalization layer, that is, for $i=1, \dots, 6$, 
		\begin{equation}
			\tilde{f}_i = \left\{
			\begin{aligned}
				&\frac{f_i - f_{\min}}{f_{\max} - f_{\min}} & \quad \text{if} \quad f_{\max} \ne f_{\min}, \\
				&1 & \quad \text{if} \quad f_{\max} = f_{\min},
			\end{aligned}
			\label{eq:3.2.3}
			\right.
		\end{equation}
		where $f_{\max} := \max\{f_i\}_{i=1}^6$ and $f_{\min} := \min\{f_i\}_{i=1}^6$. The normalization layer ensures that the magnitude of the input data does not affect the output. Following normalization, the input data proceeds through multiple hidden layers, ultimately producing a four-dimensional vector, denoted as $\bm{s}=[s_1, s_2, s_3, s_4]$. 
		In addition, we construct a hard constraint layer
		\begin{align}
			\bm{w}^{NN} = \bm{w}^* + \sum_{j=1}^4 s_j\bm{v}_j, 
			\label{eq:3.2.4}
		\end{align}
		where $\bm{w}^*$ is a fixed vector meeting these two constraints and $\bm{v}_j$ is the basis vector of the null space of the matrix $\begin{bmatrix}
			1 & 1 & 1 & 1 & 1 & 1 \\
			-5 & -3 & -1 & 1 & 3 & 5 
		\end{bmatrix}$. 
		Now we get the weights $\bm{w}^{NN}$ as the final interpolation weights for $\hat{f}_{i+\frac{1}{2}}$, which satisfy Eq. (\ref{eq:2.6a}) and (\ref{eq:2.6b}) strictly. Fig. \ref{fig:3.1} displays the architecture of WLNN.
		
		Since $\hat{f}^+_{i+\frac{1}{2}}$ and $\hat{f}^-_{i+\frac{1}{2}}$ are symmetric with respect to the cell interface, we take $\bm{w}^{NN}=\phi(f_{i+3}^-$, $f_{i+2}^-$, $ f_{i+1}^-$, $ f_{i}^-$, $ f_{i-1}^-$, $ f_{i-2}^-$; $ \theta)$ for interpolation of $\hat{f}^-_{i+\frac{1}{2}}$. Upon proper training, WLNN is expected to generate the desired weights for the 6-point finite difference scheme to solve the PDE accurately. 
		
		\subsection{Training dataset and loss function}
		\label{section3.3}
		
		Consider the scalar hyperbolic equation in Eq. (\ref{eq:2.1}) with external force $g_\lambda(x, t)$ ($\lambda$ is a parameter):  
		\begin{align}
			\frac{\partial u}{\partial t} + \frac{\partial f(u)}{\partial x} = g_\lambda(x, t). 
		\end{align}
		
		We divide the spatial and temporal intervals into $N_{x}$ and $N_{t}$ equal segments, respectively. Set $R_\lambda(x_i, t_n) := -\frac{\hat{f}_\lambda(x_{i+\frac{1}{2}}, t_n) - \hat{f}_\lambda(x_{i-\frac{1}{2}}, t_n)}{\Delta x} + g_\lambda(x_i, t_n)$. For each $\lambda$, if the exact solution $u_\lambda$ (depending on $g_\lambda$) is available, we can directly create a set of training data $\{ \{u_{\lambda, j}^n\}_{j=i-2}^{i+3}, R_{\lambda, i}^n\}$, where $u_{\lambda,j}^n:=u_{\lambda}(x_j, t^n)$ and $R_{\lambda, i}^n:= R_\lambda(x_i, t_n)$ ($1 \le i \le N_x$, $1 \le n \le N_t$); otherwise, we first need to compute the numerical solution on a finer mesh using a high-order scheme and then obtain the data on the current mesh by downsampling. Our goal is to train a neural network to generate suitable weights, which result in a good approximation of $R_\lambda(x_i, t_n)$ on the current mesh. 
		
		The loss function is given below: 
		\begin{equation}\label{eq:3.2.4a}
			\mathcal{L}(R_{\lambda, i}^{n}, R_{\lambda, i}^{n, NN}) = \frac{1}{N_{t} N_{x}} \sum_{n = 1}^{N_{t}} \sum_{i = 1}^{N_{x}} \left(R_{\lambda, i}^{n, NN} - R_{\lambda, i}^n\right)^2,
		\end{equation}
		where 
		\[
		\begin{aligned}
			& R_{\lambda, i}^{n, NN} = -\frac{\hat{f}_\lambda^{NN}(x_{i+\frac{1}{2}}, t_n) - \hat{f}_\lambda^{NN}(x_{i-\frac{1}{2}}, t_n)}{\Delta x} + g_\lambda(x_i, t_n),	\\
			& \hat{f}_\lambda^{NN}(x_{i\pm\frac{1}{2}}, t_n) = \hat{f}_{\lambda}^{NN, +}(x_{i \pm \frac{1}{2}}, t_n) + \hat{f}_{\lambda}^{NN, -}(x_{i \pm \frac{1}{2}}, t_n). 
		\end{aligned}
		\]
		Note that $\hat{f}_\lambda^{NN, \pm}(x_{i+\frac{1}{2}}, t_n)$ is computed with the weights generated by the neural network, i.e.,   
		\begin{subequations}
			\begin{align}
				& \bm{w}^{NN, +} = \phi (f_{\lambda,i-2}^{n,+}, f_{\lambda,i-1}^{n,+},  f_{\lambda,i}^{n,+}, f_{\lambda,i+1}^{n,+}, f_{\lambda,i+2}^{n,+}, f_{\lambda,i+3}^{n, +}), \\ 
				& \hat{f}_{\lambda}^{NN, +}(x_{i + \frac{1}{2}}, t_n) = \sum_{j = 1}^6 w_{j}^{NN, +} f_{\lambda,i-3+j}^{n,+},  \\
				& \bm{w}^{NN, -} = \phi (f_{\lambda,i+3}^{n,-}, f_{\lambda,i+2}^{n,-},  f_{\lambda,i+1}^{n,-}, f_{\lambda,i}^{n,-}, f_{\lambda,i-1}^{n,-}, f_{\lambda,i-2}^{n,-}), \\
				& \hat{f}_{\lambda}^{NN, -}(x_{i + \frac{1}{2}}, t_n) = \sum_{j = 1}^6 w_j^{NN, -} f_{\lambda,i+4-j}^{n,-}. 
			\end{align}
		\end{subequations}
		In practice, we need to prepare the training data $\{ \{u_{\lambda, j}^n\}_{j=i-2}^{i+3}, R_{\lambda, i}^n\}_{i,n}$ for different $\lambda$, and adopt the following loss function: 
		\begin{align}
			\mathcal{L}(R_{\lambda_s, i}^{n}, R_{\lambda_s, i}^{n, NN}) = \frac{1}{S N_{t} N_{x}} \sum_{s=1}^S \sum_{n = 1}^{N_{t}} \sum_{i = 1}^{N_{x}} \left(R_{\lambda_s, i}^{n, NN} - R_{\lambda_s, i}^n\right)^2,
			\label{eq:3.53}
		\end{align}
		where $S$ is the mini-batch size of the parameter $\lambda$. 
		
\begin{figure}
			\centering
			\includegraphics[scale=0.5]{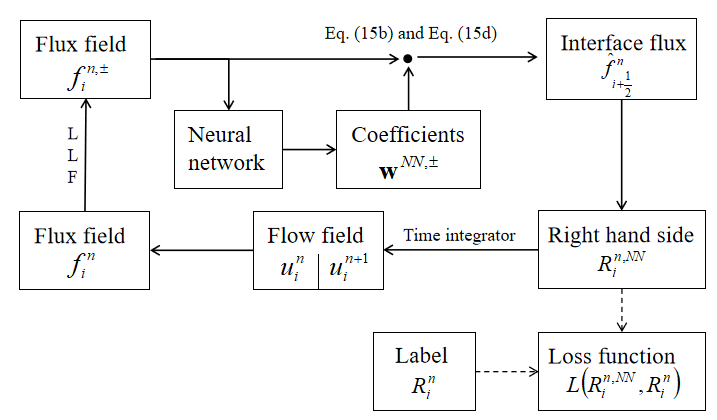}
			\caption{Learning framework and forward solver workflows.}
			\label{fig:3.2}
		\end{figure}

		Fig. \ref{fig:3.2} shows the training framework and forward solver workflow. During training, different training functions $u_i^n$ are generated, and the local flux field $ f_i^{n,\pm}$ is fed into the network to obtain the corresponding weights $\bm{w}^{NN, \pm}$. Then we get the numerical interface flux $\hat{f}_{i+1/2}^n$. By optimizing the loss function $\mathcal{L}(R_{i}^{n}, R_{i}^{n, NN})$ (the dashed part), the neural network is trained to produce the weights that result in a high-accuracy numerical scheme. During inference, we follow the solid line to apply the trained model in computation.
		
		\section{Numerical tests}
		\label{section4}
		This section presents several examples to show the performance of the WLNN method. We compare the numerical results of the WLNN method with the classical 6th-order central (CE6) scheme and the 5th-order upwind (UP5) scheme.  
		
		The third-order TVD Runge-Kutta scheme is used for time discretization\cite{ref17}:
		\begin{equation}
			\begin{aligned}
				& \bm{U}^1 = \bm{U}^n + \Delta t R(\bm{U}^n),   \\
				& \bm{U}^2 = \frac{1}{4}(3\bm{U}^n + \bm{U}^1 + \Delta t R(\bm{U}^1)),  \\
				& \bm{U}^{n+1} = \frac{1}{3}(\bm{U}^n + 2 \bm{U}^2 + 2\Delta t R(\bm{U}^2)).	\notag
			\end{aligned}
		\end{equation}
		
		The overall setting in all experiments is summarized as follows.
		\begin{itemize}
			\item \textbf{Environment} The experiments are performed in Python 3.9 with Tensorflow and use the CUDA 11.3 toolkit for GPU. 
			\item \textbf{Network setting and optimizer} The network consists of an input layer, $2$ hidden layers with nodes of $50$, and an output layer. The $\tanh$ activation function is used for the hidden layers. We train the network with the Adam optimizer \cite{ref18}. The initial learning rate is $\eta_0 = 10^{-3}$, which decays exponentially with the number of optimizations.
		\end{itemize}
		
		\subsection{Numerical examples for the 1D nonlinear scalar equation}
		\label{sec:4.1}
		The following equation is solved:
		\begin{align}
			\begin{cases}
				\frac{\partial u}{\partial t} + \frac{\partial f(u)}{\partial x} = g, \quad x\in[-2, 2],~t\in[0,1],\\
				u(x, 0) = u_0(x),  \\
				u(-2,t) = h_1(t),\quad u(2,t) = h_2(t),
				\label{eq:4.1}
			\end{cases}
		\end{align}
		where $f(u) = u^2$. To generate the dataset, we first set $u_{\lambda}(x, t) = A \exp(-(\frac{x+kt}{\sigma})^2)$ with the parameter $\lambda=\{ A, k, \sigma \}$. The external force $g(x, t) = g_{\lambda}$, boundary value $h_i(t) = h_{i,\lambda}$, $i=1, 2$ and initial value $u_0(x) = u_{0, \lambda}$ are calculated from $u_{\lambda}$. We randomly generate 50 groups of $\lambda=\{ A, k, \sigma \}$: $A \sim \mathcal{U}[0.5, 1)$, $k \sim \mathcal{U}[0, 0.5)$, $\sigma \sim \mathcal{U}[0.2, 0.3)$. And we set $N_x = 60, N_t = 40$. The number of training epochs is 1000, and the training is terminated when loss drops to $1.3\times10^{-7}$.
		
		\subsubsection{Approximate dispersion relation}
		\label{sec:4.1.1}
\begin{figure}[htbp]
			\centering
			\subfloat[]
			{\includegraphics[width=0.48\textwidth]{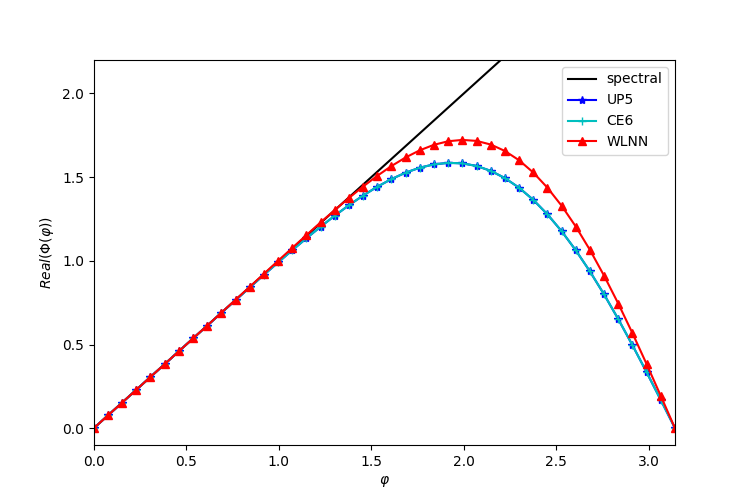}\label{fig:4.1.1a}}
			\hspace{-0.2cm}
			\subfloat[]
			{\includegraphics[width=0.48\textwidth]{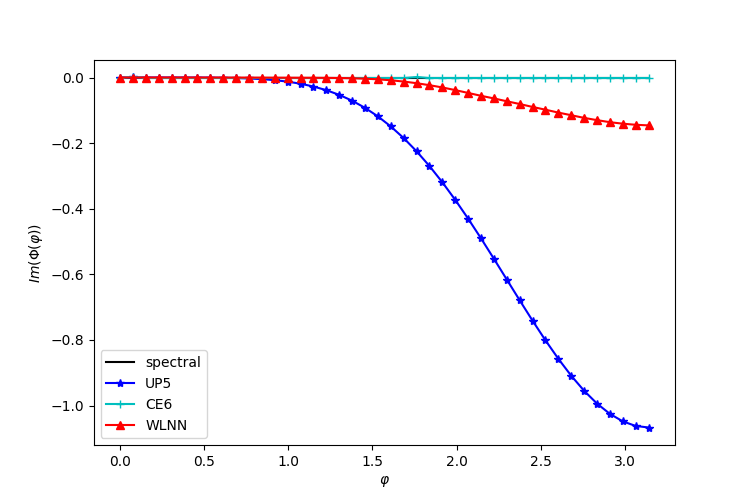}\label{fig:4.1.1b}}
			\caption{Approximate dispersion relation for the selected schemes. Real parts of $\Phi(\varphi)$ are shown in (a) and imaginary parts of $\Phi(\varphi)$ are shown in (b).}
			\label{fig:4.1.1}
		\end{figure}

		We compare the spectral properties of different schemes. Because of the nonlinearity of the WLNN method, we use the approximate dispersion relation (ADR) proposed by Pirozzoli \cite{ref19} to investigate the spectral properties, where the numerical dispersion and dissipation correspond to the deviation of the real parts and image parts of the modified wavenumber $\Phi(\varphi)$ from $\varphi$, respectively. Fig. \ref{fig:4.1.1} shows the approximated spectral properties of the CE6, UP5, and WLNN methods. The performance of all the schemes is similar for small wavenumbers. For large wavenumbers ($1.3<\varphi<2.0$), the WLNN scheme is closer to the exact spectrum in Fig. \ref{fig:4.1.1a}, which indicates it has slightly lower numerical dispersion than the UP5 and CU6 schemes. The numerical dissipation of the WLNN scheme is gradually increasing when the wavenumber $\varphi>1.6$, which is obviously smaller than the UP5 in Fig. \ref{fig:4.1.1b}. The numerical dissipation for $\varphi>2$ is not a fatal weakness and is beneficial to eliminate high-frequency oscillations \cite{ref20}.
		
		\subsubsection{Test 1-1}
		\label{sec:4.1.2}
\begin{figure}[H]
			\centering
			\includegraphics[scale=0.4]{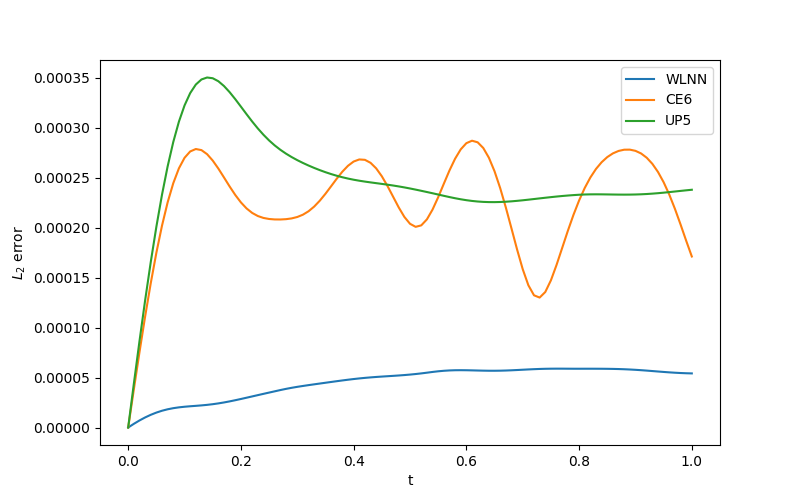}
			\caption{Test 1-1: The $L_2$ errors of WLNN, CE6 and UP5 over time.}
			\label{fig:4.1.2}
		\end{figure}

		In this test, we use the trained neural network to produce the weights of the 6-point scheme. The exact solution is $u_{\mathrm{test}}(x, t) = A \exp(-(\frac{x+kt}{\sigma})^2)$ with $\{ A, k, \sigma \} = \{0.75, 0.25, 0.25\}$. The $L_2$ errors over time are shown in Fig. \ref{fig:4.1.2}. The single-step error of WLNN is smaller than that of CE6 and UP5. In addition, we compute the mean $L_2$ errors over time, which are displayed in Table \ref{tab:4.1.1}. The mean $L_2$ error of WLNN is clearly smaller than CE6 and UP5.
		
		\begin{table}[h]
			\centering
			\begin{tabular}{cccc}
				\toprule
				Methods & WLNN & CE6 &  UP5 \\
				\midrule
				Mean $L_2$ errors & $4.5168 \times 10^{-5}$  & $2.2441 \times 10^{-4}$ & $2.4623 \times 10^{-4}$ \\
				\bottomrule
			\end{tabular}
			\caption{Test 1-1: The mean $L_2$ errors in $t \in[0, 1]$.}
			\label{tab:4.1.1}
		\end{table}
		
		\subsubsection{Test 1-2}
		\label{sec:4.1.3}
\begin{figure}[H]
			\centering
			\includegraphics[scale=0.4]{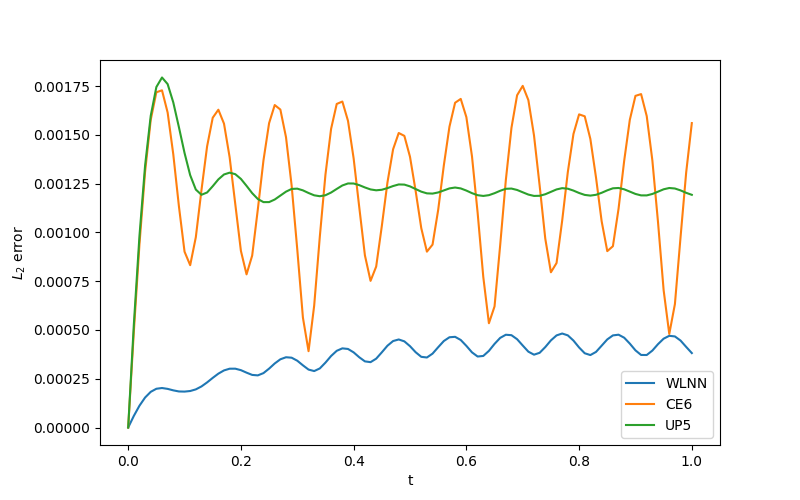}
			\caption{Test 1-2: The $L_2$ errors of WLNN, CE6 and UP5 over time.}
			\label{fig:4.1.3}
		\end{figure}

		We show that the trained neural network is also applicable to the scenario with the parameter beyond the training dataset. We set the exact solution $u_{\mathrm{test}}(x, t) = A \exp(-(\frac{x+kt}{\sigma})^2)$ with $\{A, k, \sigma\} = \{1.1, 0.7, 0.18\}$. The WLNN scheme exhibits minimal single-step error compared with CE6 and UP5. The $L_2$ errors over time are shown in Fig. \ref{fig:4.1.3}, and the average $L_2$ errors are listed in Table \ref{tab:4.1.2}. We see that WLNN has better performance than CE6 and UP5. 
		
		\begin{table}[h]
			\centering
			\begin{tabular}{cccc}
				\toprule
				Method & WLNN & CE6 & UP5 \\
				\midrule
				Mean $L_2$ errors & $3.5505 \times 10^{-4}$  & $1.2178 \times 10^{-3}$ & $1.2254 \times 10^{-3}$ \\
				\bottomrule
			\end{tabular}
			\caption{Test 1-2: The mean $L_2$ errors in $t \in[0, 1]$.}
			\label{tab:4.1.2}
		\end{table}
		
		\subsubsection{Test 1-3}
		\label{sec:4.1.4}
\begin{figure}[htbp]
			\centering
			\subfloat[]
			{\includegraphics[width=0.5\textwidth]{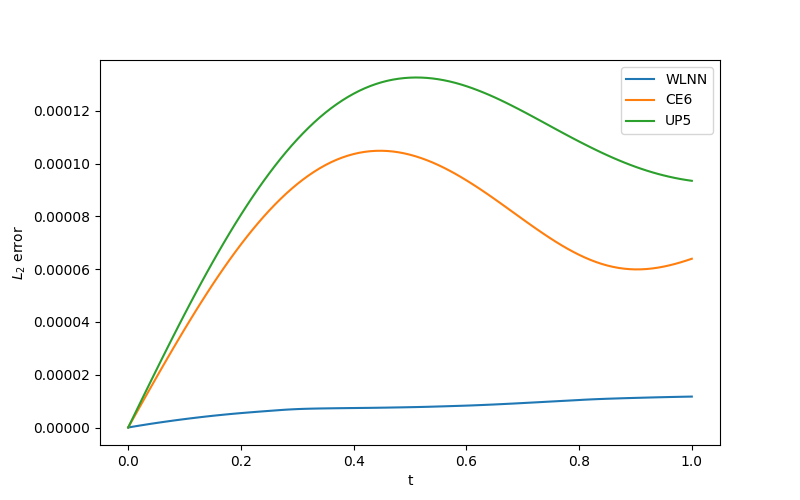}\label{fig:4.1.4a}}
			\hspace{-5mm}
			\subfloat[]
			{\includegraphics[width=0.5\textwidth]{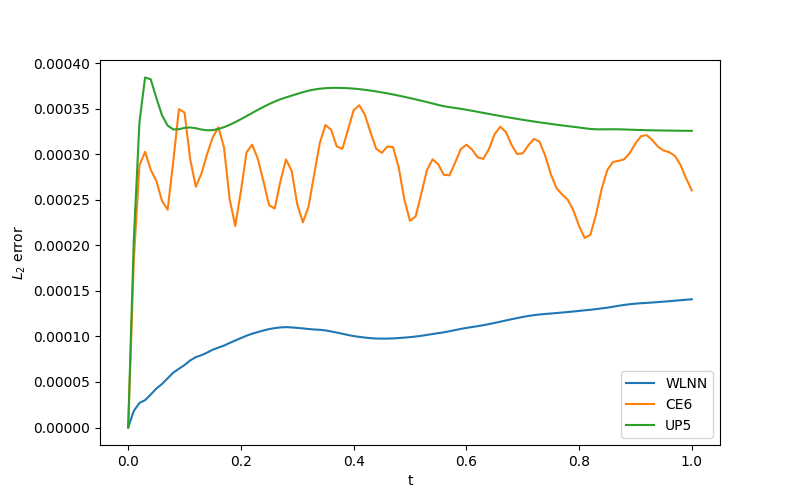}\label{fig:4.1.4b}}
			\caption{Test 1-3: The $L_2$ errors of WLNN, CE6, and UP5 over time with different $\gamma$. $\gamma = 0.1$ is shown in (a) and $\gamma = 5$ is shown in (b).}
			\label{fig:4.1.4}
		\end{figure}

		This example shows that we can also apply the trained neural network to PDEs with different coefficients. We replace the convective term in Eq. (\ref{eq:4.1}) with $f(u) = \gamma u^2$, where $\gamma$ is an adjustable parameter. As mentioned earlier, the neural network is trained for the case $\gamma = 1$, but here we conduct tests for the cases where $\gamma = 0.1$ and $\gamma = 5$, respectively. We set the exact solution $u_{\mathrm{test}}(x, t) = A \exp(-(\frac{x+kt}{\sigma})^2)$ with $\{A, k, \sigma\} = \{0.75, 0.25, 0.25\}$. The evolution of $L_2$ errors for these three schemes is depicted in Fig. \ref{fig:4.1.4}. Regardless of whether $\gamma = 0.1$ or $ \gamma = 5$, the WLNN method outperforms the other two schemes.
		
		\subsection{Numerical examples for the 3D nonlinear scalar equation}
		\label{sec:4.2}
		We consider a 3D problem on the domain $\Omega = [-2, 2] \times [-2, 2] \times [-2, 2]$: 
		\begin{align}
			\begin{cases}
				\frac{\partial u}{\partial t} + \frac{\partial f}{\partial x} + \frac{\partial f}{\partial y} + \frac{\partial f}{\partial z} = g, & \quad (x, y, z) \in \Omega, t\in\left(0, 1 \right], \\
				u(x, y, z, t) = u_b(x, y, z, t), & \quad (x, y, z) \in \partial \Omega, t\in(0, 1],	\\
				u(x, y, z, 0) = u_0(x, y, z), & \quad (x, y, z) \in \Omega,
			\end{cases}
			\label{eq:4.2}
		\end{align}
		where $f = u^2$. To generate the dataset, we set $u_{\lambda}(x, y, z, t) = A\exp(-\frac{(x+k_1t)^2 +(y+k_2t)^2 +(z+k_3t)^2}{\sigma^2})$ with $\lambda = \{A, k_1, k_2, k_3, \sigma\}$. The external force $g(x, y, z, t) = g_\lambda$, the boundary value $u_b(x, y, z, t) = u_{b, \lambda}$ and the initial value $u_0(x, y, z) = u_{0, \lambda}$ are determined by $u_{\lambda}$. We randomly choose $10$ groups of parameters $\lambda = \{A, k_1, k_2, k_3, \sigma\}$ : $A \sim \mathcal{U}[0.5, 1)$, $k_1 \sim \mathcal{U}[0, 0.5)$, $k_2 \sim \mathcal{U}[0, 0.5)$, $k_3 \sim \mathcal{U}[0, 0.5)$, $\sigma \sim \mathcal{U}[0.2, 0.3)$. 
		We set $N_x=40$ and $N_t = 20$. The training epoch is $1000$, and we stop the training when the loss is less than $1.501\times10^{-8}$.
		
		\subsubsection{Test 2-1}
		\label{sec:4.2.1}
\begin{figure}[H]
			\centering
			\includegraphics[scale=0.45]{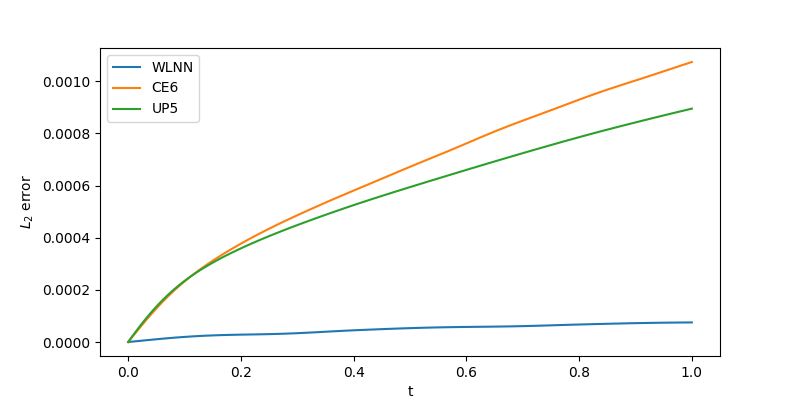}
			\caption{Test 2-1: The $L_2$ errors of WLNN, CE6 and UP5 over time.}
			\label{fig:4.2.1}
		\end{figure}

		In this test, we apply the trained neural network to solving Eq. (\ref{eq:4.2}). We select the parameters $\{A, k_1, k_2, k_3, \sigma\} = \{0.75, 0.25, 0.25, 0.25, 0.25\}$ for the exact solution. 
		The $L_2$ errors of CE6, UP5 and WLNN are displayed in Fig. \ref{fig:4.2.1}, and the mean $L_2$ errors are shown in Table \ref{tab:4.2.1}. The errors of WLNN are much smaller than those of CE6 and UP5. 
		
		\begin{table}[h]
			\centering
			\begin{tabular}{cccc}
				\toprule
				Method & WLNN & CE6 &  UP5 \\
				\midrule
				Mean $L_2$ errors & $4.7581 \times 10^{-5}$ & $6.4368 \times 10^{-4}$ & $5.6321 \times 10^{-4}$ \\
				\bottomrule
			\end{tabular}
			\caption{Test 2-1: The mean $L_2$ errors in $t \in[0, 1]$.}
			\label{tab:4.2.1}
		\end{table}
		
		\subsubsection{Test 2-2}
		\label{sec:4.2.2}
\begin{figure}[H]
			\centering
			\includegraphics[scale=0.45]{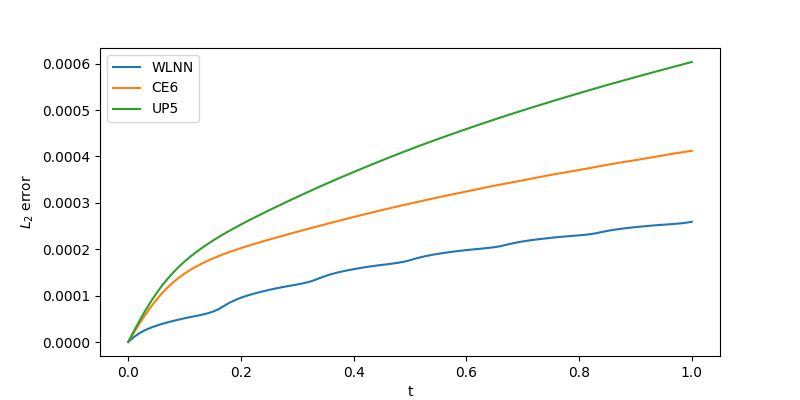}
			\caption{Test 2-2: The $L_2$ errors of WLNN, CE6 and UP5 over time.}
			\label{fig:4.2.2}
		\end{figure}

		We also assess the trained neural network's ability to address problems with parameters not included in the training dataset. For instance, we set $\{A, k_1, k_2, k_3, \sigma\} = \{1.1$, $0.6$, $0.6$, $0.6$, $0.32\}$. According to Fig. \ref{fig:4.2.2} and Table \ref{tab:4.2.2}, the $L_2$ errors of the WLNN are lower compared to CE6 and UP5.
		
		\begin{table}[h]
			\centering
			\begin{tabular}{cccc}
				\toprule
				Method & WLNN & CE6 & UP5 \\
				\midrule
				Mean $L_2$ errors & $1.6281 \times 10^{-4}$  & $2.8066 \times 10^{-4}$ & $3.8976 \times 10^{-4}$ \\
				\bottomrule
			\end{tabular}
			\caption{Test 2-2: The mean $L_2$ errors in $t \in[0, 1]$.}
			\label{tab:4.2.2}
		\end{table}
		
		\subsection{Numerical example for the 2D Euler equations}
		\label{sec:4.3}
		We apply our method to solve the 2D compressible Euler equations: 
		\begin{align}
			\frac{\bm{\partial U}}{\partial t} + \frac{\partial \bm{F}}{ \partial x} + \frac{\partial \bm{G}}{ \partial y}  = 0, \quad (x, y, t) \in \Omega \times (0, T],
		\end{align}
		where $\bm{U}=[\rho, \rho u, \rho v, E]^T$, $F = [\rho u, \rho u^2 + p, \rho u v, u(E+p)]^T$ and $G = [\rho v, \rho u v, \rho v^2 + p, v(E+p)]$. Here, $\rho$, $(u, v)$, $E$, and $p$ represent the density, velocity, total energy, and pressure, respectively.
		This system describes the conservation laws of mass, momentum, and total energy, where the total energy for an ideal polytropic gas is defined as $E = \frac{p}{(\gamma - 1)} + \frac{1}{2} \rho (u^2 + v^2)$.
		
		We consider a problem from Yee et al. \cite{ref21}, where an isentropic vortex is convected by an inviscid free stream. This numerical example is used to test the dissipative properties of different schemes. The computations are performed on the flow domain $\Omega = [0, 10] \times [0, 10]$. The initial flow field is given by
		\begin{align*}
			& \rho = \left[ 1 - \frac{(\gamma-1)\beta^2}{8\gamma \pi^2} e^{(1-r^2)} \right]^{\frac{1}{\gamma - 1}}, \quad r^2 = \Bar{x}^2 + \Bar{y}^2, \quad \Bar{x} = x - x_{vc}, \quad \Bar{y} = y - y_{vc},\\
			& (u, v) = (u_0, v_0) + \frac{\beta}{2\pi} e^{\frac{1}{2}(1-r^2)}(-\bar{y}, \bar{x}), \quad p = \rho^{\gamma},
		\end{align*}
		where $(x_{vc}, y_{vc})$ is the central coordinates of the initial vortex. The periodic boundary conditions are implemented. 
		
		We randomly choose $20$ groups of parameters $\lambda = \{x_{vc}, y_{vc}, u_0, v_0, \beta \}$ : $x_{vc} \sim \mathcal{U}[4, 6)$, $y_{vc} \sim \mathcal{U}[4, 6)$, $u_0 \sim \mathcal{U}[-1, 1)$, $v_0 \sim \mathcal{U}[-1, 1)$, $\beta \sim \mathcal{U}[2, 5)$. To generate a high accurate dataset on the mesh with $N_x = N_y = 20$, we first solve the problem on a finer grid with $N_x = N_y = 80$ and record the data $\{ \bm{U}_{i,j}, \bm{R}_{i, j} \}_{i, j=1, \dots ,80}$ at $t = \{0, 0.5, \dots, 10\}$. Averaging the data in $4\times4$ blocks, the downsampling dataset $\{ \bm{U}_{i,j}^n, \bm{R}_{i, j}^n \}_{i, j=1, \dots ,20}^{n=0, \dots, 20}$ is obtained. The number of epochs is $600$, and we stop the training when the total loss drops to $1.25 \times 10^{-5}$.
		
		\subsubsection{Test 3}
\begin{figure}
			\centering
			\includegraphics[scale=0.5]{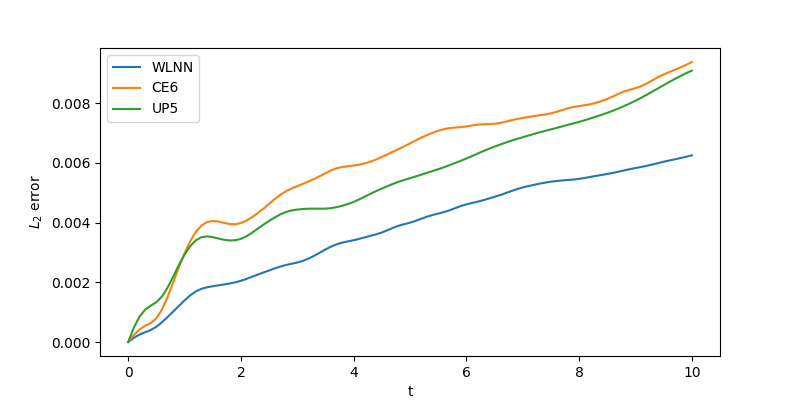}
			\caption{Test 3: The $L_2$ errors of the density $\rho(t)$ of WLNN, CE6 and UP5 over time.}
			\label{fig:4.3.1}
		\end{figure}

		We choose the exact solution with $\{x_{vc}, y_{vc}, u_0, v_0, \beta \} = \{5, 5, 1, 1, 5\}$, and carry out the simulation. The $L_2$ errors of the density $\rho(t)$ over time are shown in Fig. \ref{fig:4.3.1}. The mean $L_2$ errors in $t \in [0, 10]$ are presented in Table \ref{tab:4.3.1}. 
		We see that the errors of WLNN are smaller than CE6 and UP5.
		
		\begin{table}[h]
			\centering
			\begin{tabular}{cccc}
				\toprule
				Method & WLNN  & CE6 & UP5  \\
				\midrule
				Mean $L_2$ errors & $3.7690 \times 10^{-3}$  & $6.0425 \times 10^{-3}$  & $5.4129 \times 10^{-3}$ \\
				\bottomrule
			\end{tabular}
			\caption{Test 3: The mean $L_2$ errors of the density $\rho$ in $t \in[0, 10]$.}
			\label{tab:4.3.1}
		\end{table}
		
\begin{figure}[htbp]
			\centering
			\subfloat[\label{fig:4.3.2a}]{
				\includegraphics[scale=0.25]{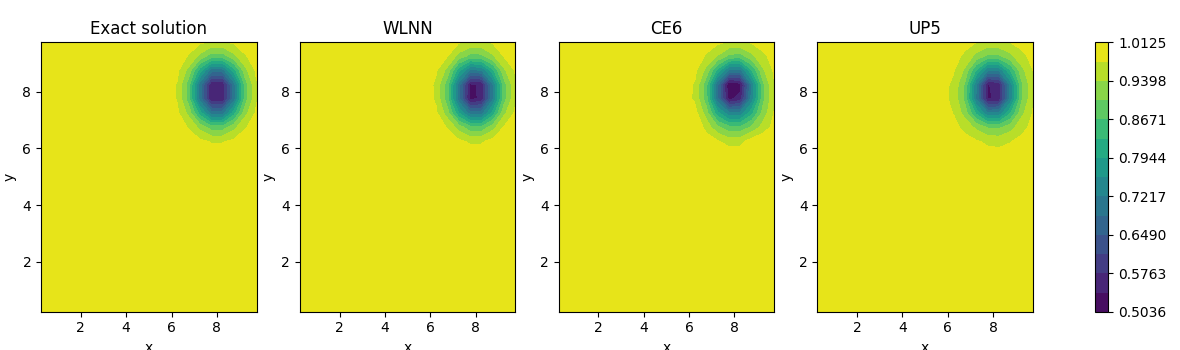}}
			\\ \vspace{-0.2cm}
			\subfloat[\label{fig:4.3.2b}]{
				\includegraphics[scale=0.25]{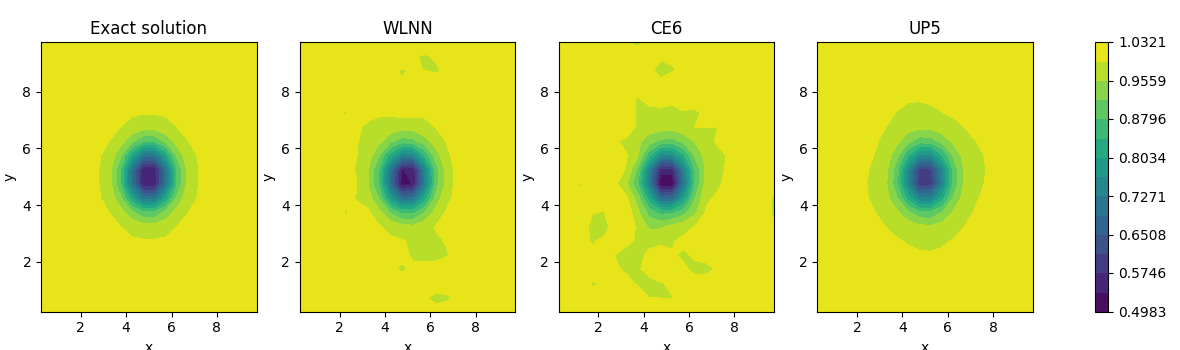}}
			\caption{Test 3: Density contours of the WLNN, CE6 and UP5 at different times. The solution is plotted with 15 contours. The density contours at $t = 3$ are shown in (a) and at $t = 10$ are shown in (b).}
			\label{fig:4.3.2}
		\end{figure}

		In Fig. \ref{fig:4.3.2}, we depict the numerical density contours computed by these three schemes in the $x-y$ plane at $t = 3 $ and $ t = 10$. At $t = 3$, the density contours of WLNN, CE6, and UP5 are similar to the exact solution. At $ t = 10$, CE6 shows deformation at the edge of the vortex, and UP5 gradually weakens the amplitude of the vortex due to the dissipative property of the scheme. The WLNN method has not encountered the aforementioned problems. It not only maintains the shape of the vortex well but also retains low dissipation.
		
		\subsection{Numerical example for the 3D Euler equations}
		\label{sec:4.4}
		We consider the 3D compressible Euler problem: 
		\begin{align}
			\frac{\bm{\partial U}}{\partial t} + \frac{\partial \bm{F}}{ \partial x} + \frac{\partial \bm{G}}{ \partial y} + \frac{\partial \bm{H}}{ \partial z} = 0, \quad (x, y, z, t) \in \Omega \times (0, T], 
		\end{align}
		with $U=[\rho, \rho u, \rho v, \rho w, E]^T$, $F = [\rho u, \rho u^2 + p, \rho v u, \rho w u, u(E+p)]^T$, $G = [\rho v, \rho u v, \rho v^2 + p, \rho w v, v(E+p)]$, $H = [\rho w, \rho u w, \rho v w, \rho w^2+p, w(E+p)]$, where $\rho$, $(u, v, w)$, and, $E$ are the density, velocity, pressure and total energy. Based on the ideal gas model, the total energy is $E = \frac{p}{(\gamma - 1)} + \frac{1}{2} \rho (u^2+v^2+w^2)$.
		
		We consider the inviscid Taylor-Green vortex, which is a benchmark case for the study of transitions to turbulence with complex manifold structures. Same as \cite{ref22}, we set the initial conditions: 
		\[
		\begin{aligned}
			& \rho = 1, \quad  u = \sin(x) \cos(y) \cos(z), \quad v = -\cos(x)\sin(y)\cos(z), \\ 
			& w = 0, \quad p = 100 + \frac{\left[ \cos(2x)+\cos(2y)\right]\left[\cos(2z)+2\right] - 2}{16},
		\end{aligned}
		\]
		where the pressure is so high that the problem is nearly incompressible. The domain is $\Omega = [0, 2\pi]^3$. The periodic boundary conditions are enforced in all directions.
		
		To generate high-accurate dataset on the mesh with $N = N_x = N_y = N_z = 32$, we solve this problem by WENO5-JS \cite{ref17} scheme together with the 3rd-order TVD Runge-Kutta scheme on a finer mesh with $N=128$ and recorder data $\{ \bm{U}_{i,j, k}, \bm{R}_{i, j, k} \}_{i, j, k=1, \dots,128}$ at $t = \{0, 0.1, \dots, 0.8\}$. The downsampling data $\{ \bm{U}_{i, j, k}^n, \bm{R}_{i, j, k}^n \}_{i, j, k=1, \dots,32}^{ n=0, \dots, 8}$ are obtained by averaging the data in $4\times4\times4$ blocks. The number of epochs is $1000$, and the total loss for training drops to $8.109 \times 10^{-6}$.
		\subsubsection{Test 4}
\begin{figure}[H]
			\centering
			\subfloat[]
			{\includegraphics[width=0.5\textwidth]{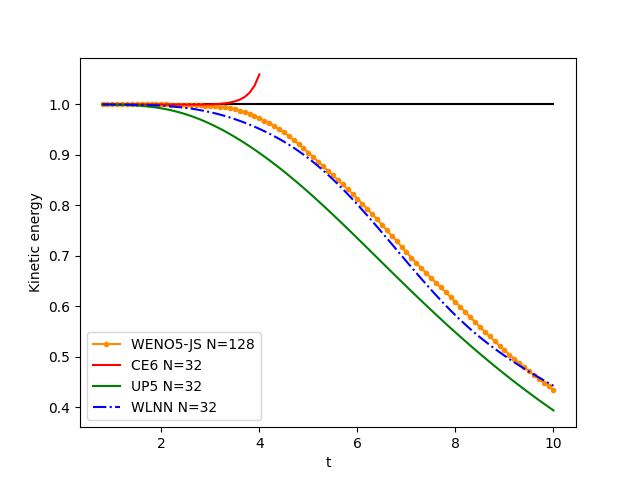}\label{fig:4.4.1a}}
			\hspace{-5mm}
			\subfloat[]
			{\includegraphics[width=0.5\textwidth]{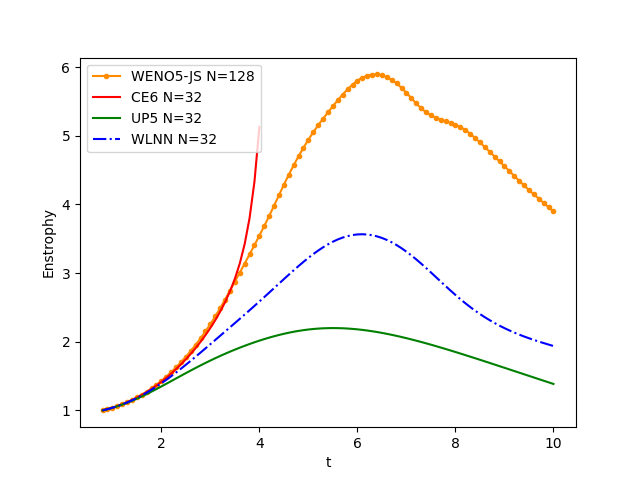}\label{fig:4.4.1b}}
			\caption{Test 4: The mean kinetic energy (a) and enstrophy (b) of WLNN, CE6 and UP5 with mesh grid $N = 32$ over time.}
			\label{fig:4.4.1}
		\end{figure}

		After training, we want to test the prediction ability of the WLNN method. The initial condition of the test problem is downsampling data at $t = 0.8$. For comparison, we conduct the simulation by using the CE6 scheme and the UP5 scheme on the same mesh with $N=32$. The result calculated by the WENO5-JS scheme with $N = 128$ is used as the reference solution.
		As with \cite{ref24}, the stability of different methods for severely under-resolved motions, as well as the conservation of kinetic energy and enstrophy growth, are concerned.
		The temporal evolution of the mean kinetic energy $K(t)$ and the enstrophy $\varOmega(t)$ is shown in Fig. \ref{fig:4.4.1}, which are normalized by $K(0.8)$ and $\varOmega(0.8)$ respectively. 
		From Fig. \ref{fig:4.4.1a}, the non-dissipative CE6 scheme maintains the conservation of kinetic energy when $t<3.7$, but the simulation is unstable after $t>4.2$, resulting in calculation failure. 
		The UP5 scheme and WLNN method preserve the kinetic energy when $t < 2$. However, the kinetic energy of these two methods gradually decays when $t > 2$ due to the numerical dissipation. The WLNN method has a slower attenuation rate than the UP5 scheme and is closer to the WENO5-JS scheme. The enstrophy of the WLNN method also performs better than the UP5 scheme in Fig. \ref{fig:4.4.1b}.
		
\begin{figure}[H]
			\centering
			\subfloat[]
			{\includegraphics[width=0.5\textwidth]{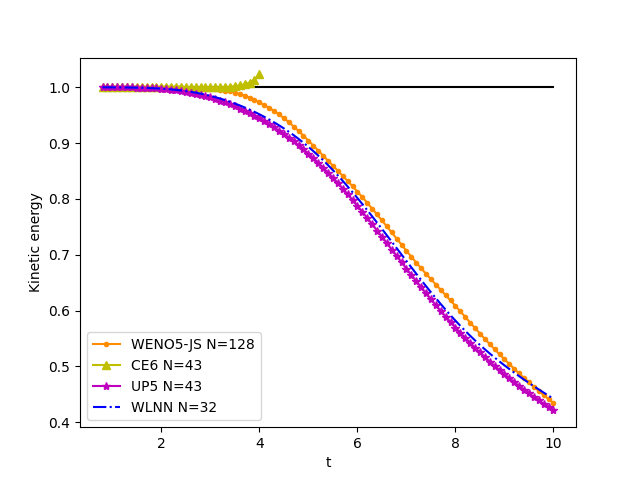}\label{fig:4.4.2a}}
			\hspace{-5mm}
			\subfloat[]
			{\includegraphics[width=0.5\textwidth]{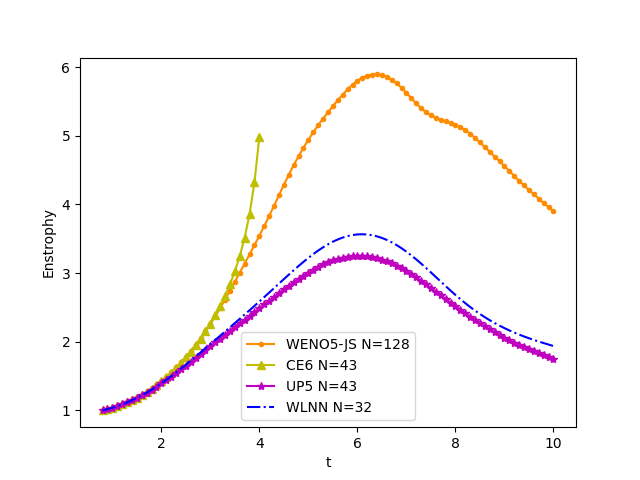}\label{fig:4.4.2b}}
			\caption{Test 4: The mean kinetic energy (a) and enstrophy (b) of WLNN, CE6 and UP5 over time.}
			\label{fig:4.4.2}
		\end{figure}

		In order to compare the efficiency, we record the running time of different schemes with $N = 32$. Since the CE6 scheme quickly becomes unstable, we only compare the WLNN method with the UP5 scheme. The running times of the UP5 scheme and WLNN method are $528.2$s and $1057.1$s, respectively. 
		Even though the WLNN method takes nearly twice as much time as the UP5 scheme for $N = 32$, the accuracy of WLNN is noticeably higher.
		We also evaluate the results of the UP5 scheme of $N = 43$ shown in Fig. \ref{fig:4.4.2}, which exhibits similar performances with the WLNN method of $N = 32$.
		The running time of the UP5 scheme with $N=43$ is $2014.5$s, which is nearly twice as long as the WLNN of $N=32$.

		\section{Conclusion}
		\label{section5}
		In this paper, we propose high-accuracy numerical schemes learned by the neural network to compute the interface flux, which can automatically adjust the weights (or equivalently, the reconstruction coefficients) according to the local flux information. We apply the WLNN method to the nonlinear scalar hyperbolic problems and the 2D/3D Euler equations. The experimental results show our method is accurate and can be generalized to the scenarios with parameters beyond the training dataset and the PDEs with different coefficients as well. The errors of the WLNN scheme are smaller than those of the CE6 and UP5 schemes. In the simulation of the inviscid Taylor-Green vortex, the WLNN runs faster than UP5 under similar accuracy, while the CE6 is unstable. We only consider the problem with smooth solutions, and the extension to the problem with non-smooth solutions will be envisioned in future work.  
		
		\section*{Acknowledgments}
		Y. G. was partially supported by the NSFC Grant G0592370101. 
		H. S. was partially supported by the NSFC Grant 62231016.
		L. X. was supported by the NSFC Grants 12431015 and 62231016. 
		J. Z. was supported by the NSFC Grant 12201024.

		\bibliographystyle{elsarticle-num-names}
		\bibliography{references}

\begin{thebibliography}{34}
\expandafter\ifx\csname natexlab\endcsname\relax\def\natexlab#1{#1}\fi
\providecommand{\url}[1]{\texttt{#1}}
\providecommand{\href}[2]{#2}
\providecommand{\path}[1]{#1}
\providecommand{\DOIprefix}{doi:}
\providecommand{\ArXivprefix}{arXiv:}
\providecommand{\URLprefix}{URL: }
\providecommand{\Pubmedprefix}{pmid:}
\providecommand{\doi}[1]{\href{http://dx.doi.org/#1}{\path{#1}}}
\providecommand{\Pubmed}[1]{\href{pmid:#1}{\path{#1}}}
\providecommand{\bibinfo}[2]{#2}
\ifx\xfnm\relax \def\xfnm[#1]{\unskip,\space#1}\fi
\bibitem[{Giraldo and Restelli(2008)}]{ref00}
\bibinfo{author}{F.~X. Giraldo}, \bibinfo{author}{M.~Restelli},
\newblock \bibinfo{title}{A study of spectral element and discontinuous
  galerkin methods for the {Navier--Stokes} equations in nonhydrostatic
  mesoscale atmospheric modeling: Equation sets and test cases},
\newblock \bibinfo{journal}{Journal of Computational Physics}
  \bibinfo{volume}{227} (\bibinfo{year}{2008}) \bibinfo{pages}{3849--3877}.
\bibitem[{Yi(2018)}]{ref01}
\bibinfo{author}{T.-H. Yi},
\newblock \bibinfo{title}{Time integration of unsteady nonhydrostatic equations
  with dual time stepping and multigrid methods},
\newblock \bibinfo{journal}{Journal of Computational Physics}
  \bibinfo{volume}{374} (\bibinfo{year}{2018}) \bibinfo{pages}{873--892}.
\bibitem[{Suresh and Huynh(1997)}]{ref1}
\bibinfo{author}{A.~Suresh}, \bibinfo{author}{H.~T. Huynh},
\newblock \bibinfo{title}{Accurate monotonicity-preserving schemes with
  {Runge--Kutta} time stepping},
\newblock \bibinfo{journal}{Journal of Computational Physics}
  \bibinfo{volume}{136} (\bibinfo{year}{1997}) \bibinfo{pages}{83--99}.
\bibitem[{Fang et~al.(2013)Fang, Li, and Lu}]{ref2}
\bibinfo{author}{J.~Fang}, \bibinfo{author}{Z.~Li}, \bibinfo{author}{L.~Lu},
\newblock \bibinfo{title}{An optimized low-dissipation monotonicity-preserving
  scheme for numerical simulations of high-speed turbulent flows},
\newblock \bibinfo{journal}{Journal of Scientific Computing}
  \bibinfo{volume}{56} (\bibinfo{year}{2013}) \bibinfo{pages}{67--95}.
\bibitem[{Li et~al.(2013)Li, Leng, and He}]{ref3}
\bibinfo{author}{X.-l. Li}, \bibinfo{author}{Y.~Leng}, \bibinfo{author}{Z.-w.
  He},
\newblock \bibinfo{title}{Optimized sixth-order monotonicity-preserving scheme
  by nonlinear spectral analysis},
\newblock \bibinfo{journal}{International Journal for Numerical Methods in
  Fluids} \bibinfo{volume}{73} (\bibinfo{year}{2013})
  \bibinfo{pages}{560--577}.
\bibitem[{Jin et~al.(2018)Jin, Liao, and Cai}]{ref6}
\bibinfo{author}{Y.~Jin}, \bibinfo{author}{F.~Liao}, \bibinfo{author}{J.~Cai},
\newblock \bibinfo{title}{Optimized low-dissipation and low-dispersion schemes
  for compressible flows},
\newblock \bibinfo{journal}{Journal of Computational Physics}
  \bibinfo{volume}{371} (\bibinfo{year}{2018}) \bibinfo{pages}{820--849}.
\bibitem[{Liang and Fu(2024)}]{ref7}
\bibinfo{author}{T.~Liang}, \bibinfo{author}{L.~Fu},
\newblock \bibinfo{title}{A new type of non-polynomial based {TENO} scheme for
  hyperbolic conservation laws},
\newblock \bibinfo{journal}{Journal of Computational Physics}
  \bibinfo{volume}{497} (\bibinfo{year}{2024}) \bibinfo{pages}{112618}.
\bibitem[{Pinkus(1999)}]{ref51}
\bibinfo{author}{A.~Pinkus},
\newblock \bibinfo{title}{Approximation theory of the {MLP} model in neural
  networks},
\newblock \bibinfo{journal}{Acta Numerica} \bibinfo{volume}{8}
  (\bibinfo{year}{1999}) \bibinfo{pages}{143--195}.
\bibitem[{Leshno et~al.(1993)Leshno, Lin, Pinkus, and Schocken}]{ref52}
\bibinfo{author}{M.~Leshno}, \bibinfo{author}{V.~Y. Lin},
  \bibinfo{author}{A.~Pinkus}, \bibinfo{author}{S.~Schocken},
\newblock \bibinfo{title}{Multilayer feedforward networks with a nonpolynomial
  activation function can approximate any function},
\newblock \bibinfo{journal}{Neural networks} \bibinfo{volume}{6}
  (\bibinfo{year}{1993}) \bibinfo{pages}{861--867}.
\bibitem[{Raissi et~al.(2019)Raissi, Perdikaris, and Karniadakis}]{ref53}
\bibinfo{author}{M.~Raissi}, \bibinfo{author}{P.~Perdikaris},
  \bibinfo{author}{G.~E. Karniadakis},
\newblock \bibinfo{title}{Physics-informed neural networks: A deep learning
  framework for solving forward and inverse problems involving nonlinear
  partial differential equations},
\newblock \bibinfo{journal}{Journal of Computational Physics}
  \bibinfo{volume}{378} (\bibinfo{year}{2019}) \bibinfo{pages}{686--707}.
\bibitem[{Yu et~al.(2018)}]{ref54}
\bibinfo{author}{B.~Yu}, et~al.,
\newblock \bibinfo{title}{The deep {Ritz} method: a deep learning-based
  numerical algorithm for solving variational problems},
\newblock \bibinfo{journal}{Communications in Mathematics and Statistics}
  \bibinfo{volume}{6} (\bibinfo{year}{2018}) \bibinfo{pages}{1--12}.
\bibitem[{Cai et~al.(2021)Cai, Chen, and Liu}]{ref55}
\bibinfo{author}{Z.~Cai}, \bibinfo{author}{J.~Chen}, \bibinfo{author}{M.~Liu},
\newblock \bibinfo{title}{Least-squares {ReLU} neural network ({LSNN}) method
  for linear advection-reaction equation},
\newblock \bibinfo{journal}{Journal of Computational Physics}
  \bibinfo{volume}{443} (\bibinfo{year}{2021}) \bibinfo{pages}{110514}.
\bibitem[{Li et~al.(2024{\natexlab{a}})Li, Jiang, Sun, Xu, and Zhou}]{ref531}
\bibinfo{author}{H.~Li}, \bibinfo{author}{S.~Jiang}, \bibinfo{author}{W.~Sun},
  \bibinfo{author}{L.~Xu}, \bibinfo{author}{G.~Zhou},
\newblock \bibinfo{title}{A model-data asymptotic-preserving neural network
  method based on micro-macro decomposition for gray radiative transfer
  equations},
\newblock \bibinfo{journal}{Communications in Computational Physics}
  \bibinfo{volume}{35} (\bibinfo{year}{2024}{\natexlab{a}})
  \bibinfo{pages}{1155--1193}.
\bibitem[{Li et~al.(2024{\natexlab{b}})Li, Jiang, Sun, Xu, and Zhou}]{ref532}
\bibinfo{author}{H.~Li}, \bibinfo{author}{S.~Jiang}, \bibinfo{author}{W.~Sun},
  \bibinfo{author}{L.~Xu}, \bibinfo{author}{G.~Zhou},
\newblock \bibinfo{title}{Macroscopic auxiliary asymptotic preserving neural
  networks for the linear radiative transfer equations},
\newblock \bibinfo{journal}{arXiv preprint arXiv:2403.01820}
  (\bibinfo{year}{2024}{\natexlab{b}}).
\bibitem[{Sirignano and Spiliopoulos(2018)}]{ref56}
\bibinfo{author}{J.~Sirignano}, \bibinfo{author}{K.~Spiliopoulos},
\newblock \bibinfo{title}{{DGM}: A deep learning algorithm for solving partial
  differential equations},
\newblock \bibinfo{journal}{Journal of Computational Physics}
  \bibinfo{volume}{375} (\bibinfo{year}{2018}) \bibinfo{pages}{1339--1364}.
\bibitem[{Patel et~al.(2022)Patel, Manickam, Trask, Wood, Lee, Tomas, and
  Cyr}]{ref57}
\bibinfo{author}{R.~G. Patel}, \bibinfo{author}{I.~Manickam},
  \bibinfo{author}{N.~A. Trask}, \bibinfo{author}{M.~A. Wood},
  \bibinfo{author}{M.~Lee}, \bibinfo{author}{I.~Tomas}, \bibinfo{author}{E.~C.
  Cyr},
\newblock \bibinfo{title}{Thermodynamically consistent physics-informed neural
  networks for hyperbolic systems},
\newblock \bibinfo{journal}{Journal of Computational Physics}
  \bibinfo{volume}{449} (\bibinfo{year}{2022}) \bibinfo{pages}{110754}.
\bibitem[{Bar-Sinai et~al.(2019)Bar-Sinai, Hoyer, Hickey, and Brenner}]{ref13}
\bibinfo{author}{Y.~Bar-Sinai}, \bibinfo{author}{S.~Hoyer},
  \bibinfo{author}{J.~Hickey}, \bibinfo{author}{M.~P. Brenner},
\newblock \bibinfo{title}{Learning data-driven discretizations for partial
  differential equations},
\newblock \bibinfo{journal}{Proceedings of the National Academy of Sciences}
  \bibinfo{volume}{116} (\bibinfo{year}{2019}) \bibinfo{pages}{15344--15349}.
\bibitem[{Bezgin et~al.(2021)Bezgin, Schmidt, and Adams}]{ref12}
\bibinfo{author}{D.~A. Bezgin}, \bibinfo{author}{S.~J. Schmidt},
  \bibinfo{author}{N.~A. Adams},
\newblock \bibinfo{title}{A data-driven physics-informed finite-volume scheme
  for nonclassical undercompressive shocks},
\newblock \bibinfo{journal}{Journal of Computational Physics}
  \bibinfo{volume}{437} (\bibinfo{year}{2021}) \bibinfo{pages}{110324}.
\bibitem[{Kochkov et~al.(2021)Kochkov, Smith, Alieva, Wang, Brenner, and
  Hoyer}]{ref15}
\bibinfo{author}{D.~Kochkov}, \bibinfo{author}{J.~A. Smith},
  \bibinfo{author}{A.~Alieva}, \bibinfo{author}{Q.~Wang},
  \bibinfo{author}{M.~P. Brenner}, \bibinfo{author}{S.~Hoyer},
\newblock \bibinfo{title}{Machine learning--accelerated computational fluid
  dynamics},
\newblock \bibinfo{journal}{Proceedings of the National Academy of Sciences}
  \bibinfo{volume}{118} (\bibinfo{year}{2021}) \bibinfo{pages}{e2101784118}.
\bibitem[{Tompson et~al.(2017)Tompson, Schlachter, Sprechmann, and
  Perlin}]{ref8}
\bibinfo{author}{J.~Tompson}, \bibinfo{author}{K.~Schlachter},
  \bibinfo{author}{P.~Sprechmann}, \bibinfo{author}{K.~Perlin},
\newblock \bibinfo{title}{Accelerating eulerian fluid simulation with
  convolutional networks},
\newblock in: \bibinfo{booktitle}{International conference on machine
  learning}, \bibinfo{organization}{PMLR}, \bibinfo{year}{2017}, pp.
  \bibinfo{pages}{3424--3433}.
\bibitem[{Wang et~al.(2024)Wang, Cao, Huang, Liu, Hu, Luo, Song, Zhao, Liu, Sun
  et~al.}]{ref81}
\bibinfo{author}{H.~Wang}, \bibinfo{author}{Y.~Cao},
  \bibinfo{author}{Z.~Huang}, \bibinfo{author}{Y.~Liu},
  \bibinfo{author}{P.~Hu}, \bibinfo{author}{X.~Luo}, \bibinfo{author}{Z.~Song},
  \bibinfo{author}{W.~Zhao}, \bibinfo{author}{J.~Liu},
  \bibinfo{author}{J.~Sun}, et~al.,
\newblock \bibinfo{title}{Recent advances on machine learning for computational
  fluid dynamics: A survey},
\newblock \bibinfo{journal}{arXiv preprint arXiv:2408.12171}
  (\bibinfo{year}{2024}).
\bibitem[{Zhuang et~al.(2021)Zhuang, Kochkov, Bar-Sinai, Brenner, and
  Hoyer}]{ref14}
\bibinfo{author}{J.~Zhuang}, \bibinfo{author}{D.~Kochkov},
  \bibinfo{author}{Y.~Bar-Sinai}, \bibinfo{author}{M.~P. Brenner},
  \bibinfo{author}{S.~Hoyer},
\newblock \bibinfo{title}{Learned discretizations for passive scalar advection
  in a two-dimensional turbulent flow},
\newblock \bibinfo{journal}{Physical Review Fluids} \bibinfo{volume}{6}
  (\bibinfo{year}{2021}) \bibinfo{pages}{064605}.
\bibitem[{Bezgin et~al.(2022)Bezgin, Schmidt, and Adams}]{ref10}
\bibinfo{author}{D.~A. Bezgin}, \bibinfo{author}{S.~J. Schmidt},
  \bibinfo{author}{N.~A. Adams},
\newblock \bibinfo{title}{{WENO3-NN}: A maximum-order three-point data-driven
  weighted essentially non-oscillatory scheme},
\newblock \bibinfo{journal}{Journal of Computational Physics}
  \bibinfo{volume}{452} (\bibinfo{year}{2022}) \bibinfo{pages}{110920}.
\bibitem[{Xue et~al.(2022)Xue, Xia, Li, and Yuan}]{ref11}
\bibinfo{author}{Z.~Xue}, \bibinfo{author}{Y.~Xia}, \bibinfo{author}{C.~Li},
  \bibinfo{author}{X.~Yuan},
\newblock \bibinfo{title}{A simplified multilayer perceptron detector for the
  hybrid {WENO} scheme},
\newblock \bibinfo{journal}{Computers \& Fluids} \bibinfo{volume}{244}
  (\bibinfo{year}{2022}) \bibinfo{pages}{105584}.
\bibitem[{de~Lara and Ferrer(2022)}]{ref9}
\bibinfo{author}{F.~M. de~Lara}, \bibinfo{author}{E.~Ferrer},
\newblock \bibinfo{title}{Accelerating high order discontinuous {Galerkin}
  solvers using neural networks: 1d {Burgers'} equation},
\newblock \bibinfo{journal}{Computers \& Fluids} \bibinfo{volume}{235}
  (\bibinfo{year}{2022}) \bibinfo{pages}{105274}.
\bibitem[{de~Lara and Ferrer(2023)}]{ref91}
\bibinfo{author}{F.~M. de~Lara}, \bibinfo{author}{E.~Ferrer},
\newblock \bibinfo{title}{Accelerating high order discontinuous {Galerkin}
  solvers using neural networks: {3D} compressible {Navier-Stokes} equations},
\newblock \bibinfo{journal}{Journal of Computational Physics}
  \bibinfo{volume}{489} (\bibinfo{year}{2023}) \bibinfo{pages}{112253}.
\bibitem[{Toro(2013)}]{ref16}
\bibinfo{author}{E.~F. Toro}, \bibinfo{title}{Riemann solvers and numerical
  methods for fluid dynamics: a practical introduction},
  \bibinfo{publisher}{Springer Science \& Business Media},
  \bibinfo{year}{2013}.
\bibitem[{Jiang and Shu(1996)}]{ref17}
\bibinfo{author}{G.-S. Jiang}, \bibinfo{author}{C.-W. Shu},
\newblock \bibinfo{title}{Efficient implementation of weighted {ENO} schemes},
\newblock \bibinfo{journal}{Journal of Computational Physics}
  \bibinfo{volume}{126} (\bibinfo{year}{1996}) \bibinfo{pages}{202--228}.
\bibitem[{Kingma(2014)}]{ref18}
\bibinfo{author}{D.~P. Kingma},
\newblock \bibinfo{title}{Adam: A method for stochastic optimization},
\newblock \bibinfo{journal}{arXiv preprint arXiv:1412.6980}
  (\bibinfo{year}{2014}).
\bibitem[{Pirozzoli(2006)}]{ref19}
\bibinfo{author}{S.~Pirozzoli},
\newblock \bibinfo{title}{On the spectral properties of shock-capturing
  schemes},
\newblock \bibinfo{journal}{Journal of Computational Physics}
  \bibinfo{volume}{219} (\bibinfo{year}{2006}) \bibinfo{pages}{489--497}.
\bibitem[{Shen(2023)}]{ref20}
\bibinfo{author}{H.~Shen},
\newblock \bibinfo{title}{A class of {ENO} schemes with adaptive order for
  solving hyperbolic conservation laws},
\newblock \bibinfo{journal}{Computers \& Fluids} \bibinfo{volume}{266}
  (\bibinfo{year}{2023}) \bibinfo{pages}{106050}.
\bibitem[{Yee et~al.(1999)Yee, Sandham, and Djomehri}]{ref21}
\bibinfo{author}{H.~C. Yee}, \bibinfo{author}{N.~D. Sandham},
  \bibinfo{author}{M.~J. Djomehri},
\newblock \bibinfo{title}{Low-dissipative high-order shock-capturing methods
  using characteristic-based filters},
\newblock \bibinfo{journal}{Journal of Computational Physics}
  \bibinfo{volume}{150} (\bibinfo{year}{1999}) \bibinfo{pages}{199--238}.
\bibitem[{Deng et~al.(2020)Deng, Jiang, Xiao, and Yan}]{ref22}
\bibinfo{author}{X.~Deng}, \bibinfo{author}{Z.-H. Jiang},
  \bibinfo{author}{F.~Xiao}, \bibinfo{author}{C.~Yan},
\newblock \bibinfo{title}{Implicit large eddy simulation of compressible
  turbulence flow with {PnTm- BVD} scheme},
\newblock \bibinfo{journal}{Applied Mathematical Modelling}
  \bibinfo{volume}{77} (\bibinfo{year}{2020}) \bibinfo{pages}{17--31}.
\bibitem[{Johnsen et~al.(2010)Johnsen, Larsson, Bhagatwala, Cabot, Moin, Olson,
  Rawat, Shankar, Sj{\"o}green, Yee et~al.}]{ref24}
\bibinfo{author}{E.~Johnsen}, \bibinfo{author}{J.~Larsson},
  \bibinfo{author}{A.~V. Bhagatwala}, \bibinfo{author}{W.~H. Cabot},
  \bibinfo{author}{P.~Moin}, \bibinfo{author}{B.~J. Olson},
  \bibinfo{author}{P.~S. Rawat}, \bibinfo{author}{S.~K. Shankar},
  \bibinfo{author}{B.~Sj{\"o}green}, \bibinfo{author}{H.~C. Yee}, et~al.,
\newblock \bibinfo{title}{Assessment of high-resolution methods for numerical
  simulations of compressible turbulence with shock waves},
\newblock \bibinfo{journal}{Journal of Computational Physics}
  \bibinfo{volume}{229} (\bibinfo{year}{2010}) \bibinfo{pages}{1213--1237}.

\end{thebibliography}
		
	\end{document}